\newcommand{\text}[1]{\quad\mbox{#1}\quad}
\def\pont{\hspace{-6pt}{\bf.\ }}
\def\beq{\begin{equation}}\def\eeq{\end{equation}}
\def\eps{\varepsilon}
\def\qed{\ifhmode\unskip\nobreak\fi\quad\ifmmode\Box\else$\Box$\fi}

\documentclass[12pt]{article} \textheight 8in\textwidth
6in \oddsidemargin 0in\evensidemargin 0in
\usepackage{graphicx}
\usepackage{amsfonts}
\usepackage{epsf}

\newtheorem{theorem}{Theorem}
\newtheorem{lemma}[theorem]{Lemma}

\begin{document}
\title{Ramsey number of a connected triangle matching}
\author{Andr\'as Gy\'arf\'as\thanks{Research supported in part by
the OTKA Grant No. K104343.}\\\\[-0.8ex]
\small Alfr\'ed R\'enyi Institute of Mathematics\\[-0.8ex]
\small Hungarian Academy of Sciences\\[-0.8ex]
\small Budapest, P.O. Box 127\\[-0.8ex]
\small Budapest, Hungary, H-1364 \\\small
\texttt{gyarfas.andras@renyi.mta.hu} \and G\'{a}bor N. S\'ark\"ozy
\thanks{Research supported in part by
the National Science Foundation under Grant No. DMS-0968699 and by OTKA Grant No. K104343.}\\[-0.8ex]
\small Alfr\'ed R\'enyi Institute of Mathematics\\[-0.8ex]
\small Hungarian Academy of Sciences\\[-0.8ex]
\small Budapest, P.O. Box 127\\[-0.8ex]
\small Budapest, Hungary, H-1364\\
and\\
\small Computer Science Department\\[-0.8ex]
\small Worcester Polytechnic Institute\\[-0.8ex]
\small Worcester, MA, USA 01609\\[-0.8ex]
\small \texttt{gsarkozy@cs.wpi.edu}}

\maketitle

\begin{abstract}

We determine the $2$-color Ramsey number of a {\em connected}
triangle matching $c(nK_3)$ which is any connected graph containing
$n$ vertex disjoint triangles. We obtain that
$R(c(nK_3),c(nK_3))=7n-2$, somewhat larger than in the classical
result of Burr, Erd\H os and Spencer for a triangle matching,
$R(nK_3,nK_3)=5n$. The motivation is to determine the Ramsey number
$R(C_n^2,C_n^2)$ of the square of a cycle $C_n^2$. We apply our
Ramsey result for connected triangle matchings to show that the
Ramsey number of an ``almost" square of a cycle $C_n^{2,c}$ (a cycle
of length $n$ in which all but at most a constant number $c$ of
short diagonals are present) is asymptotic to $7n/3$.

\end{abstract}

\section{Introduction}

Denote by $\delta(G)$ the minimum degree in a graph $G$. $K_n$ is
the complete graph on $n$ vertices and $K_{n,n}$ is the complete
bipartite graph between two sets of $n$ vertices each. If $G_1, G_2,
\ldots , G_r$ are graphs, then the Ramsey number $R(G_1, G_2, \ldots
, G_r)$ is the smallest positive integer $n$ such that in any
edge-coloring with colors $1,2,\ldots,r$, for some $i$ the edges of
color $i$ contain a subgraph isomorphic to $G_i$. In this paper we
will deal with 2- and 3-color Ramsey numbers (so $r=2$ or 3) and we
will think of color 1 as red, color 2 as blue and color 3 (if it
exists) as white.

Among well-known early results in generalized Ramsey theory is the
exact value of $R(nK_2,nK_2)=3n-1$ determined by Cockayne and
Lorimer {\cite{CL}, and $R(nK_3,nK_3)=5n$ determined by Burr, Erd\H
os and Spencer \cite{BES}. Here $nG$ denotes $n$ vertex disjoint
copies of the graph $G$. It turned out in many applications that it
is important to study the case when $nK_2$, the matching,  is
replaced by a {\em connected matching}, $c(nK_2)$, defined as any
connected graph containing $nK_2$ (see for example applications
\cite{BBGGYS,FL,GYLSS,rlogr,GRSSz1,GRSSz3,GYS,GYSSZdiam,L}). The
2-color Ramsey number of connected matchings is the same as the
Ramsey number of matchings (in fact, \cite{GGY} proves more) and one
of the key arguments of \cite{GRSSz1} was that this remains true for
3 colors as well. However, for more than $3$ colors the Ramsey
numbers of matchings and connected matchings are different. For
example, $R(nK_2,nK_2,nK_2,nK_2)=5n-3$ \cite{CL}, but
$R(c(nK_2),c(nK_2),c(nK_2),c(nK_2))>6n-3$ when $2n-1$ is divisible
by three. This can be seen by the $4$-coloring obtained from the
parallel classes of an affine plane of order $3$ by replacing each
point with a point set of size ${2n-1\over 3}$.

In this paper we look at the connected version of the ``matching of
triangles''. Let $c(nK_3)$ denote any connected graph containing $n$
vertex disjoint triangles. We shall prove that here already the
$2$-color Ramsey number of $c(nK_3)$ is different from its counterpart $nK_3$.

\begin{theorem}\pont \label{conntr} For $n\ge 2$, $R(c(nK_3),c(nK_3))=7n-2$.
\end{theorem}
(While we have $R(nK_3,nK_3)=5n$.) To see that
$R(c(nK_3),c(nK_3))>7n-3$, consider pairwise disjoint sets $A,B,C$
such that $V(K_{7n-3})=A\cup B \cup C$, $|A|=|B|=3n-1$ and
$|C|=n-1$. Edges inside $A$ and inside $B$ are red, all other edges
are blue. In this coloring there is no monochromatic $c(nK_3)$ (in
fact there is not even an $nK_3$ in blue).

To prove that $R(c(nK_3),c(nK_3))\le 7n-2$, we need the Ramsey number of connected
triangle matchings versus ordinary matchings that might be interesting
on its own.

\begin{lemma} \pont \label{tr} For $1\le m \le n$, $R(c(nK_3),mK_2)=3n+m-1$.
\end{lemma}

Notice that $R(c(nK_3),mK_2)>3n+m-2$ is shown by the disjoint sets
$A,X$ such that $V(K_{3n+m-2})=A\cup X$, $|A|=3n-1,|X|=m-1$ and
edges inside $A$ are colored red, other edges are colored with blue.
In this coloring there is no red $nK_3$ or blue $mK_2$.

The motivation of Theorem \ref{conntr} comes from the effort to
determine or estimate the 2-color Ramsey number $R(C_n^2,C_n^2)$ where of $C_n^2$ is the
square of the cycle on $n$ vertices, i.e. the cycle $C_n$ with all
short diagonals (diagonals between vertices at a distance 2 on the
cycle). A very recent paper of Allen, Brightwell and Skokan \cite{ABS} gives lower bound $3n-4$ and upper bound ${20n\over 3}+o(n)$
for that Ramsey number.

Density questions for the square of a cycle  have also received a
lot of attention (for example the well-known P\'osa-Seymour problem,
see \cite{CDK}, \cite{KSS3}, \cite{LSSz}). Note also that there has
been a lot of research on the Ramsey numbers of constant maximum
degree graphs (such as $C_n^2$, where the maximum degree is 4). It
is known that for a graph $G$ on $n$ vertices with maximum degree
$\Delta$ the Ramsey number is linear, it is at most $c(\Delta)n$
\cite{CRST}. The current best bound $c(\Delta)\leq
2^{c\Delta\log{\Delta}}$ is due to Conlon, Fox and Sudakov
\cite{CFS}.

To determine $R(C_n^2,C_n^2)$ exactly for every $n$ is hopeless
since for $n=5,6$, $C_5^2=K_5, C_6^2=K_6-3K_2$ and these Ramsey
numbers are both unknown. Also, the constant $c(4)$ in the linear
bounds $c(4)n$ is very large.  However, combining Theorem
\ref{conntr} with the Regularity Lemma, we shall make a step forward
and prove the following. Let $C_n^{2,c}$ denote an ``almost" square
of the cycle $C_n$, a cycle of length $n$ in which all but at most a
constant number $c$ of short diagonals are present.

\begin{theorem}\pont \label{7/3}
For every fixed $\eta > 0$ there is a $c=c(\eta)$ so that for any
$n\ge c$ we have $R(C_n^{2,c}, C_n^{2,c}) \leq (1+\eta) 7n/3$.
\end{theorem}
It is worth noting that  Theorem \ref{7/3}, although asymptotically
sharp (shown by a similar construction as in Theorem \ref{conntr}),
does not give the right asymptotics for $R(C_n^2,C_n^2)$, where we
insist on {\em all} short diagonals. Indeed, $R(C_n^2,C_n^2)\ge
3n-4$ is proved in \cite{ABS}. Thus perhaps surprisingly removing
these constant number of diagonals makes a big difference in the
Ramsey number.

The following easy lemma from \cite{GYS} will be used. It extends
(when $\delta(G)=|V(G)|-1$) the well-known remark of Erd\H os and
Rado that in a $2$-colored complete graph there is a monochromatic
spanning tree.

\begin{lemma}\pont\label{conn}(Lemma 1.5 in \cite{GYS}) Suppose that the edges of a graph $G$ with $\delta(G)\ge {3|V(G)|\over4}$ are
$2$-colored. Then there is a monochromatic connected subgraph with
order larger than $\delta(G)$. This estimate is sharp.
\end{lemma}

Theorem \ref{conntr} and Lemma \ref{tr} are proved in Section
\ref{ramsey}. Their perturbed versions are worked out in Section
\ref{perturbedramsey}. Section \ref{reg} outlines the (rather
standard) argument how to obtain Theorem \ref{7/3} from the
Regularity and Blow-up lemmas.

\section{The proof of Lemma \ref{tr} and Theorem \ref{conntr}}\label{ramsey}

\noindent {\bf Proof of Lemma \ref{tr}. }We prove by induction on $m$. Since for $m=1$ the statement is
trivially true for every $n$, suppose we have a blue matching
$M=(m-1)K_2$ in a $2$-coloring of a $K_{3n+m-1}$ with vertex set
$V$. If there is no blue $mK_2$ then every edge $e_i\in M$ has a
vertex $p_i$ adjacent in red to all but at most one vertex of
$X=V-V(M)$ (otherwise $e_i$ could be replaced by two independent
blue edges). Also, $X$ induces a red complete graph. Since
$$|X|=3n+m-1-2(m-1)=3n-m+1\ge 2m+1,$$ we can select greedily $m-1$ pairwise
disjoint red triangles with one vertex as $p_i$ and two vertices
from $X$. Then we find red triangles greedily in the remainder of
$X$. We are guaranteed to find $n$ red triangles this way. These
triangles can be certainly included into a connected red subgraph so
we have the required $c(nK_3)$. \qed

\bigskip
\noindent {\bf Proof of Theorem \ref{conntr}. } Consider a
$2$-coloring of the edges of $K=K_{7n-2}$ with red and blue, assume
w.l.o.g. that the blue color class has only one connected component
(since one color class is connected by using Lemma \ref{conn} with
$\delta(G)=|V(G)|-1$). Since $R(nK_3,nK_3)=5n$ for $n\geq 2$, we
have a monochromatic $nK_3$. If it is blue, we are done, therefore
it is red and the red color class must define at least two connected
components. Suppose that the red components have vertex sets
$X_1,X_2,\dots,X_s$, where $s\ge 2$ and $|X_1|\ge |X_2|\ge\dots\ge
|X_s|\ge 1$. We may suppose that $|X_1|\le 5n$, otherwise we have
the required monochromatic $c(nK_3)$ from $R(nK_3,nK_3)=5n$.

\noindent {\bf Case 1.} $|X_1|\ge |X_2|\ge 3n$, $|X_1|=3n+k_1$,
$|X_2|=3n+k_2$. Since by Lemma \ref{tr} we have
$R(c(nK_3),(k_i+1)K_2)=3n+k_i$ for $i=1,2$, $X_i$ contains either a
red $c(nK_3)$ or a blue $(k_i+1)K_2$ and we are done if the first
possibility appears. Thus we have blue matchings $M_i$ of size
$k_i+1$ in $X_i$ for $i=1,2$. We can take $k_1+1$ vertices in
$X_2\setminus M_2$ and $k_2+1$ vertices in $X_1\setminus M_1$ to
form a blue $T=(k_1+k_2+2)K_3$ using the blue edges between $X_1$
and $X_2$ since from $k_1+k_2\le n-2$ it follows that
$$3n+k_1-2(k_1+1)=3n-k_1-2\ge k_2+1,$$
$$3n+k_2-2(k_2+1)=3n-k_2-2\ge k_1+1.$$
If $k_1+k_2+2 = n$ (i.e. $s=2$) we have the required blue $c(nK_3)$.
Otherwise we have $l=n-(k_1+k_2)-2$ ($>0$) vertices in $A=V-(X_1\cup
X_2)$ and we can form $l$ vertex disjoint blue triangles taking one
vertex from each of the sets $A,X_1-T,X_2-T$. We have enough
vertices for that, because
$$|X_1-T|=3n+k_1-2(k_1+1)-(k_2+1)=3n-(k_1+k_2)-3\ge n-(k_1+k_2)-2$$
and the same is true for $|X_2-T|$. Thus we have a connected blue
triangle matching of size at least $k_1+k_2+2+n-(k_1+k_2)-2=n$, as
desired finishing Case 1.

\noindent {\bf Case 2.} $|X_1|\ge 3n$, $2n\le |X_2|<3n$,
$|X_1|=3n+k_1$, $|X_2|=3n-k_2$, $1\le k_2\le n$. Again since by
Lemma \ref{tr} we have $R(c(nK_3),(k_1+1)K_2)=3n+k_1$, we may
suppose that we have a blue $M_1=(k_1+1)K_2$ in $X_1$. We transform
$M_1$ to a blue triangle matching $T=(k_1+1)K_3$ using $k_1+1$
vertices from $X_2$ and then extend $T$ using  $q=n-k_1+k_2-2$
vertices in $A=V-(X_1\cup X_2)$ and $q$ vertices from both
$X_1,X_2$. We have enough room for that because
$$3n+k_1-2(k_1+1)\ge n-k_1+k_2-2, \mbox{  } 3n-k_2-(k_1+1)\ge n-k_1+k_2-2$$
are both true since $k_2\le n$. Thus we have a connected blue
triangle matching of size at least $(k_1+1)+(n-k_1+k_2-2)=n+k_2-1\ge
n$.

\noindent {\bf Case 3.} $|X_1|\ge 3n$, $n\le |X_2|<2n$,
$|X_1|=3n+k_1$, $|X_2|=n+k_2$, $0\le k_2<n$. Again since by Lemma
\ref{tr} we have $R(c(nK_3),(k_1+1)K_2)=3n+k_1$, we may suppose that
we have a blue $M_1=(k_1+1)K_2$ in $X_1$. Furthermore, we may
suppose that $k_1< n$ otherwise $M_1$ can be transformed to a blue
triangle matching $T=(k_1+1)K_3$ using vertices from $X_2$. Since
$k_2<n$ also holds, $|V-(X_1\cup X_2)|=3n-(k_1+k_2)-2\ge n$. Thus we
have at least $n$ vertices in all of the three sets
$X_1,X_2,V-(X_1\cup X_2)$ implying that we have a connected blue
$c(nK_3)$.

\noindent {\bf Case 4.} $|X_1|\ge 3n$, $|X_2|<n$ or $|X_1|< 3n$.

From Lemma \ref{tr} we may assume $|X_1|\le 4n-2$. Indeed, otherwise
since by Lemma \ref{tr} we have $R(c(nK_3),nK_2)=4n-1$, we may
suppose that we have a blue $M_1=nK_2$ in $X_1$. This blue $M_1$ can
be transformed into a blue $c(nK_3)$ using $n$ vertices in $V-X_1$
($|X_1|<5n$ ensures that there are $n$ vertices).

Define the set $S_1$ so that $|S_1|=n$ and, starting with $X_1$, all
vertices of $X_i$ are selected before taking vertices from
$X_{i+1}$. Then, starting from the next $X_i$, define $S_2$ in the
same way. Now set
$$A=\{\cup X_i: (S_1\cup S_2) \cap X_i= \emptyset\}.$$
Then we have the following claim. $$|A|\geq n \;(\mbox{or
equivalently}\;  |V\setminus A|\leq 6n-2).$$ Indeed, this is true
either because $|X_1|\le 4n-2$ and $|X_2|<n$ so the last $X_i$ which
intersects $S_2$ satisfies $|X_i|<n$ or because $3n>|X_1|\ge |X_2|$.
But then we can select $S_3\subset A$ with $|S_3|= n$, and the blue
complete tripartite graph $[S_1,S_2,S_3]$ defines the required blue
$c(nK_3)$. \qed

\section{Perturbed version of Theorem \ref{conntr}}\label{perturbedramsey}

As in many applications of the Regularity Lemma, one has to handle a
few irregular pairs and the corresponding edges will not be present
in the reduced graph. We say that the graph $G$ on $n$ vertices is
$\varepsilon$-{\it perturbed} if it is almost complete, at most
$\varepsilon {n\choose 2}$ edges are missing. We cannot apply
Theorem \ref{conntr} in the reduced graph because in Theorem
\ref{conntr} we have a 2-colored complete graph, yet the reduced
graph will be a 2-colored $\varepsilon$-perturbed graph. Thus we
need perturbed versions of Theorem \ref{conntr} and first Lemma
\ref{tr}. It will be convenient to think of the missing edges as
edges in a third color class (white or ``invisible''), so
we move up from 2-color Ramsey numbers to 3-color Ramsey numbers.
$K_{1,t}$ denotes the star with $t$ leaves.

\begin{lemma} \pont \label{pert-tr} For $1\le m\le n$, $0\leq t\leq n$, $R(c(nK_3),mK_2,K_{1,t})\leq 3n+m-1+2t$.
\end{lemma}

\noindent{\bf Proof. } We prove by induction on $m$ as in the
non-perturbed case. The starting case, $m=1$ follows easily from a
well-known result of Corradi and Hajnal \cite{CH} (or it could also
be proved directly by an easy induction on $n$). Indeed, if there is
no blue edge, we have a red graph on $N=3n+2t$ vertices with minimum
degree at least $3n+t>{2N\over 3}$ and it contains at least $\lfloor
N/3 \rfloor \ge n$ vertex disjoint red triangles. Since our red
graph is automatically connected from the minimum degree condition,
we have the required red $c(nK_3)$. Thus, we may select a blue
matching $M=(m-1)K_2$ in a $2$-coloring of a $K_{3n+m-1+2t}$ with
vertex set $V$. We may assume that from every vertex fewer than $t$
edges are missing (or white edges). If there is no blue $mK_2$ then
every edge $e_i\in M$ again has a vertex $p_i$ that is adjacent in
blue to at most one vertex in $X=V-V(M)$. However, now $p_i$ is not
necessarily adjacent in red to all other vertices in $X$ since some
edges might be missing. But all the edges that are actually {\em
present} are indeed red to the other vertices. Furthermore, in $X$
all edges that are present are red as well. Since
$$|X|=3n+m-1+2t-2(m-1)=3n-m+1+2t\ge 2m+1+2t,$$ again we can select greedily $m-1$ pairwise
disjoint red triangles with one vertex as $p_i$ and two vertices
from $X$. Indeed, $p_i$ is still adjacent in red to more than
$(2t+3)-t=t+3>t$ vertices in $X$ but then there is a (red) edge
among these neighbors, giving a red triangle as desired. Then we
find red triangles greedily in the remainder of $X$ similarly.
Finally, we find the $n$-th red triangle in the remainder of $X$ as
follows. Select an arbitrary remaining vertex of $X$. Since it has
more than $t$ neighbors left in $X$, there is an edge among these
neighbors and all edges are red in $X$. The red graph spanned by $X$
is connected because $|X|/2 > t$, thus the $n$ red triangles form a
$c(nK_3)$.\qed

We will also need a perturbed version of the classical result of
Burr, Erd\H os and Spencer, $R(nK_3,nK_3)=5n$.

\begin{lemma} \pont \label{pert-burr} For $n\ge 2$, $1\le t \le n, R(nK_3,nK_3,K_{1,t})\leq 6n-2+5t$.
\end{lemma}
\noindent{\bf Proof. } Consider the largest blue triangle matching,
remove it, then consider the largest red triangle matching in the
remainder and remove it. We have at least $5t+4$ vertices left and
there are no more monochromatic triangles. However, consider an
arbitrary vertex, it is still adjacent to at least $4t+4$ vertices
in the leftover. Then in one of the colors (say blue) it is adjacent
to at least $2t+2$ vertices. These neighbors will induce a triangle
which must be red (otherwise we get a blue triangle), a
contradiction. Indeed, consider again an arbitrary vertex from these
at least $2t+2$ vertices, it is still adjacent to at least $t+2$
vertices from these at least $2t+2$ vertices. But then there must be
an edge within these at least $t+2$ vertices, giving a triangle.

We note that more is true ($5n+ct$) but for our purposes this weaker
statement is sufficient. \qed

Next we are ready to give the perturbed version of Theorem
\ref{conntr}.

\begin{theorem}\pont \label{pert-conntr} For $n\ge 2$, $0\leq t\leq 2n/3$, $R(c(nK_3),c(nK_3),K_{1,t})\leq 7n-2+7t$.
\end{theorem}

\noindent {\bf Proof. } Again suppose we have a $2$-coloring of a
$K_{7n-2+7t}$ with vertex set $V$. We may assume that from every
vertex fewer than $t$ edges are missing (edges in the third color).
Applying Lemma \ref{conn}, we get a monochromatic (say blue)
connected component $X$ of size at least $(7n-2+7t)-t=7n-2+6t$. By
Lemma \ref{pert-burr}, since $7n-2+6t\geq 6n-2+5t$, we have a
monochromatic $nK_3$ in $X$. If it is blue, we are done, therefore
it is red and thus the red color class must define at least two
connected components within $X$. Suppose that the red components of
$V$ have vertex sets $X_1,X_2,\dots,X_s$, where $s\ge 2$ and
$|X_1|\ge |X_2|\ge\dots\ge |X_s|\ge 1$. We may suppose that
$|X_1|\le 6n-2+5t$ otherwise we have the required monochromatic
$c(nK_3)$ from Lemma \ref{pert-burr}.

\noindent {\bf Case 1.} $|X_1|\ge |X_2|\ge 3n+2t$, $|X_1|=3n+2t+k_1$,
$|X_2|=3n+2t+k_2$. Here we apply Lemma \ref{pert-tr} to the two
subgraphs induced by $X_1$ and $X_2$. We find either a red $c(nK_3)$
or a blue $(k_i+1)K_2$ in them and we are done if the first
possibility appears. Thus we have blue matchings $M_1,M_2$ of size
$k_1+1,k_2+1$, respectively.

We extend $M_1$ to a blue $(k_1+1)K_3$ by taking $k_1+1$ vertices
in $X_2-M_2$. This can be done if $|X_2-M_2|\ge k_1+1+2t-2$, extending the edges of $M_1$ one by one to blue triangles, at each step we have at most $2t-2$ vertices in $X_2-M_2$ not adjacent (in blue) to at least one of the ends of the edge to be extended. Indeed,
$$|X_2-M_2|=3n+2t+k_2-2(k_2+1)=3n+2t-k_2-2\ge k_1+1+2t-2$$
i.e $3n\ge k_1+k_2+1$ which is true since $k_1+k_2+1\le n+3t-1\le n+2n-1$ from the assumption $t\le 2n/3$. The same argument allows to extend $M_2$ to a blue $(k_2+1)K_3$ with $k_2+1$ vertices of $X-M_1$. Thus we have a blue $T=(k_1+k_2+2)K_3$ and noticing that the blue graph spanned by $X_1\cup X_2$ is connected (by $|X_1|,|X_2|>2t$, any two vertices of $X_1$ and of $X_2$ has a common blue neighbor) we are done if $k_1+k_2+2 \ge n$.

Otherwise we have $l=n+3t-(k_1+k_2)-2$
($>0$) vertices in $A=V-(X_1\cup X_2)$ and we plan to extend $T$ to $nK_3$ with $n-(k_1+k_2+2)$ vertex
disjoint blue triangles taking one vertex from each of the sets
$A,X_1-T,X_2-T$. Since $T$ is already connected, the extension will be automatically connected as well. We have enough vertices for that if all the three sets have size at least $n-(k_1+k_2+2)+2t$ (in fact two of them can be only at least $n-(k_1+k_2+2)+t$).  In our case the condition holds for $A$ with $t$ to save and for $X_i-T$ it holds with about $2n$ to save:
$$|X_1-T|=3n+2t+k_1-2(k_1+1)-(k_2+1)=$$
$$=3n+2t-(k_1+k_2)-3\ge n-(k_1+k_2+2)+2t$$
and the same is true for
$|X_2-T|$.  Thus we have a blue $c(nK_3)$.

\noindent {\bf Case 2.} $|X_1|\ge 3n+2t$, $2n+2t\le |X_2|<3n+2t$,
$|X_1|=3n+2t+k_1$, $|X_2|=3n+2t-k_2$, $1\le k_2\le n$.

\noindent Here we apply Lemma
\ref{pert-tr} to the subgraph induced by $X_1$. We may suppose we have a
blue $M_1=(k_1+1)K_2$ in $X_1$, we transform $M_1$ to a blue
triangle matching $T=(k_1+1)K_3$ using $k_1+1$ vertices from $X_2$
and then extend $T$ using  $q=n-k_1+k_2-2$ vertices in $A=V-(X_1\cup
X_2)$ and $q$ vertices from both $X_1,X_2$. We have enough room for
that because $|A|=n+3t-k_1+k_2-2\ge q+2t$ and
$$3n+2t+k_1-2(k_1+1)\ge q+2t, \mbox{  } 3n+2t-k_2-(k_1+1)\ge q+2t$$
are both true since $k_2\le n$. Thus we have a connected blue
triangle matching of size at least $(k_1+1)+q=n+k_2-1\ge
n$.

\noindent {\bf Case 3.} $|X_1|\ge 3n+2t$, $n+2t\le |X_2|<2n+2t$,
$|X_1|=3n+2t+k_1$, $|X_2|=n+2t+k_2$, $0\le k_2<n$.

\noindent Again we apply Lemma \ref{pert-tr} to the subgraph induced by $X_1$ and
select the blue $M_1=(k_1+1)K_2$ in $X_1$. We may suppose that $k_1<
n$ otherwise $M_1$ can be transformed to a blue triangle matching
$T=(k_1+1)K_3$ using vertices from $X_2$. Since $k_2<n$ also holds,
$|V-(X_1\cup X_2)|=3n+3t-(k_1+k_2)-2\ge n+2t$. Thus we have at least $n+2t$
vertices in all of the three sets $X_1,X_2,V-(X_1\cup X_2)$ implying
that we have a connected blue $c(nK_3)$.

\noindent {\bf Case 4.} $|X_1|\ge 3n+2t$, $|X_2|<n+2t$ or $|X_1|< 3n+2t$.

\noindent We may assume $|X_1|\le 4n-2+2t$, otherwise we can apply Lemma \ref{pert-tr} with $m=n$ to $X_1$ to find a blue $nK_2$ and, since from  Lemma \ref{pert-conntr} $|X_1|\le 6n-2+5t$ we have at least $7n-2+7t-(6n-2+5t)=n+2t$ vertices in $V-X_1$, the blue $nK_2$  can be
transformed into a blue $c(nK_3)$ using $n$ vertices of $V-X_1$.

If $|X_1|\ge n+2t$, take an $(n+2t)$-vertex subset $S_1\subset X_1$ then
take an $(n+2t)$-vertex set $S_2$ from $\cup_{i>1} X_i$ so that in $S_2$
we use all vertices of $X_i$ before taking vertices from $X_{i+1}$.
Define
$$A=\{\cup_{i>1}X_i : S_2 \cap X_i\ne \emptyset\}$$
Then $|X_1\cup A|\le 6n-2+6t$ either because $|X_1|\le 4n-2+2t$ and
$|X_2|<n+2t$ so the last $X_i$ which intersects $S_2$ satisfies
$|X_i|<n+2t$ or because $3n+2t>|X_1|\ge |X_2|$. Thus we can select
$S_3\subset V-A$ with $|S_3|\ge n+t$, and the blue tripartite
graph $[S_1,S_2,S_3]$ has lower bounds $n+2t,n+2t,n+t$ for its vertex classes which allows to pick the vertices of the required blue $c(nK_3)$.

If $|X_1|<n+2t$, define $S_1$ so that $|S_1|=n+t$ and all vertices of
$X_i$ are selected before taking vertices from $X_{i+1}$. Then,
starting from the next $X_i$, define $S_2$ in the same way. Now set
$$B=\{\cup X_i: (S_1\cup S_2) \cap X_i\ne \emptyset\}$$ and observe
that $|B|\le 4n+6t$ thus we can select $S_3\subset V-B$ with
$|S_3|\ge 7n-2+7t-(4n+6t)=3n-2+t\ge n+2t$, and the blue complete tripartite graph $[S_1,S_2,S_3]$ has lower bounds $n+t,n+t,n+2t$ for its vertex classes which allows to pick the vertices of the required blue $c(nK_3)$.  \qed

\section{Proof of Theorem \ref{7/3}; applying the Regularity
Lemma}\label{reg}

Next we show how to prove Theorem~\ref{7/3} from Theorem
\ref{pert-conntr}, the Regularity Lemma~\cite{Sz} and the Blow-up
Lemma. The material of this section is fairly standard by now (see
\cite{BBGGYS,GYLSS,rlogr,GRSSz1,GRSSz3,GYS,GYSSZdiam} for similar
techniques) so we omit some of the details. In particular in
\cite{GRSSz1} Section 2 follows a similar outline.

Let $e(X,Y)$ denote the number of edges between $X$ and $Y$ in a
graph $G$. For disjoint $X,Y$, we define the density
$$d(X,Y)={e(X,Y)\over|X|\cdot|Y|}\,.$$
For two disjoint subsets $A,B$ of $V(G)$, the bipartite graph with
vertex set $A\cup B$ which has all the edges of $G$ with one
endpoint in $A$ and the other in $B$ is called the pair $(A,B)$.

A pair $(A,B)$ is $\eps$-regular if for every $X\subset A$ and
$Y\subset B$ satisfying
$$|X|>\eps|A|\text{and}|Y|>\eps|B|$$
we have $$|d(X,Y)-d(A,B)|<\eps.$$ A pair $(A,B)$ is
$(\eps,\delta)$-super-regular if it is $\eps$-regular and
furthermore,$$deg(a)\geq\delta|B|\text{for all}a\in A,$$
$$\text{and}deg(b)\geq\delta|A|\text{for all}b\in B.$$

We need a $2$-edge-colored version of the Szemer\'edi Regularity
Lemma.\footnote{For background, this variant and other variants of
the Regularity Lemma see \cite{KS}.}

\begin{lemma}\label{2-reg}\pont
For every integer $m_0$ and positive $\varepsilon$, there is an
$M_0=M_0(\eps, m_0)$ such that for $n\geq M_0$ the following holds.
For any $n$-vertex  graph $G$, where $G=G_1\cup G_2$ with
$V(G_1)=V(G_2)=V$, there is a partition of $V$ into $\ell+1$
clusters $V_0,V_1,\dots,V_\ell$ such that
\begin{itemize}
\item $m_0\leq \ell\leq M_0$, $|V_1|=|V_2|=\dots=|V_\ell|$,
$|V_0|<\eps n$,
 \item apart from at most $\eps {\ell
\choose 2}$ exceptional pairs, all pairs $G_s|_{V_i\times V_j}$ are
$\eps$-regular, where $1\leq i<j\leq \ell$ and $1\leq s\leq 2$.
\end{itemize}
\end{lemma}

Our other main tool is the Blow-up Lemma (see \cite{KSSz2, KSSz3}).
It basically says that super-regular pairs behave like complete
bipartite graphs from the point of view of bounded degree subgraphs.

\begin{lemma}\label{blow}\pont
Given a graph $R$ of order $r$ and positive parameters
$\delta,\Delta$, there exists an $\eps>0$ such that the following
holds. Let $m$ be an arbitrary positive integer, and let us replace
the vertices of $R$ with pairwise disjoint $m$-sets
$V_1,V_2,\ldots,V_r$ (blowing up). We construct two graphs on the
same vertex-set $V=\cup V_i$. The graph $R(m)$ is obtained by
replacing all edges of $R$ with copies of the complete bipartite
graph $K_{m,m}$, and a sparser graph $G$ is constructed by replacing
the edges of $R$ with some $(\eps,\delta)$-super-regular pairs. If a
graph $H$ with $\Delta(H)\leq\Delta$ is embeddable into $R(m)$ then
it is already embeddable into $G$.
\end{lemma}

Actually we will need the following consequence of the Blow-up Lemma
(where $R$ is a triangle).

\begin{lemma}\pont\label{blow-up}
For every $\delta>0$ there exist an $\eps > 0$ and $m_0$ such that
the following holds. Let $G$ be a tripartite graph with tripartition
$V(G)=V_1\cup V_2\cup V_3$ such that $|V_1|=|V_2|=|V_3|=m\geq m_0$,
and let all the 3 pairs $(V_1,V_2)$, $(V_1,V_3)$, $(V_2,V_3)$ be $(\eps,
\delta)$-super-regular. Then for every pair of vertices $v_1\in V_1,
v_2\in V_3$ and for every integer $p$, $4\leq p \leq 3m$, $G$
contains an ``almost" $P_p^2$, the square of a path with $p$
vertices connecting $v_1$ and $v_2$ from which  at most two
short diagonals are missing.
\end{lemma}

We emphasize that Lemma \ref{blow-up} is true for any value of $p$
between $4$ and $3m$, not just for the ones that are divisible by
$3$. The price we pay is that two short diagonals might be missing
which is allowed in our application. Note also that an easier
approximate version of this lemma would suffice as well, but for
simplicity we use this lemma.

\noindent {\bf Proof. } We think of $G$ as having the orientation
$(V_1,V_2,V_3)$. Because of the Blow-up Lemma it is sufficient to
check the statement for the complete tripartite graph (using $r=3$
and $\Delta=4$ in the Blow-up Lemma). We distinguish three cases
depending on $p$. If $p=3k$ for some integer $2\leq k\leq m$, then
we just go around $(V_1,V_2,V_3)$ $k$ times starting with $v_1$ and
ending with $v_2$, so in this case actually no short diagonal is
missing. If $p=3k+1$ with $k\geq 1$, then we go around
$(V_1,V_2,V_3)$ $(k-1)$ times starting with $v_1$, but then in the
last round we ``turn around", i.e. we finish with vertices
$u_1,u_2,u_3,v_2$ chosen from $V_1,V_2,V_1, V_3$, respectively. Then
the only short diagonal missing is between $u_1$ and $u_3$. Finally,
if $p=3k+2$ with $k\geq 1$, then we go around $(V_1,V_2,V_3)$
$(k-1)$ times starting with $v_1$, but then in the last round we
``double up", i.e. we finish with vertices $u_1,u_2,u_3,u_4,v_2$
chosen from $V_1,V_2,V_1,V_2,V_3$, respectively. Then the only two
short diagonals missing are between $u_1$ and $u_3$ and between
$u_2$ and $u_4$. \qed

With these preparations now we are ready to prove Theorem~\ref{7/3}
from Theorem \ref{pert-conntr}. Let \beq\label{para}\eps\ll \eta \ll
1,\eeq $m_0$ sufficiently large compared to $1/\varepsilon$ and
$M_0$ obtained from Lemma~\ref{2-reg}. Suppose we have a
$2$-coloring of a complete graph with vertex set $V$,
$|V|=(1+\eta)7n/3$ (for simplicity assume that this is a
sufficiently large integer). We apply Lemma~\ref{2-reg}. We obtain a
partition of $V$, that is $V=\cup_{0\leq i\leq \ell}V_i$. We define
the following {\em reduced graph} $G^R$: The vertices of $G^R$ are
$p_1, \ldots , p_\ell$, and there is an edge between vertices $p_i$
and $p_j$ if the pair $(V_i, V_j)$ is $\eps$-regular in both colors.
The edge $p_ip_j$ is colored with the majority color in $K(V_i,
V_j)$. Thus $G^R$ is a $(1-\eps)$-dense 2-colored graph on $\ell$
vertices. Then we ``trim" $G^R$ in the standard way: there is a
subgraph $H^{R}$ on at least $(1-\sqrt{\eps})\ell$ vertices where
the maximum degree of the complement is less than $\sqrt{\eps}\ell$
(see for example Lemma 9 in \cite{GRSSz1}). In other words the third
color class does not contain a star $K_{1,t}$ with
$t=\sqrt{\eps}\ell$, as we need in Theorem \ref{pert-conntr}.

Applying Theorem \ref{pert-conntr} to $H^R$ with
$t=\sqrt{\eps}\ell$, we can get a large monochromatic (say red)
connected triangle matching in $H^R$ (and thus in $G^R$). For a
triangle $T_i, 1\leq i \leq \ell_1$ in this connected triangle
matching denote the corresponding clusters by $(V_1^i,V_2^i,V_3^i)$.
Thus (using (\ref{para})) we may assume that the number of vertices
in the union of these clusters is between $(1+\frac{\eta}{2})n$ and
$(1+\eta)n$. Next, first using the fact that this is a connected
triangle matching we find red connecting paths $P_i^R$ in $G^R$
between $T_i$ and $T_{i+1}$, $1\leq i \leq \ell_1$ (where
$T_{\ell_1+1}=T_1$) and then from these connecting paths $P_i^R$ we
can find vertex disjoint red connecting paths $P_i$ in the original
graph between a vertex of $V_3^i$ to a vertex of $V_1^{i+1}$. This
procedure is rather standard by now, see for example \cite{GRSSz1}
for the details.

These connecting paths $P_i$ will be part of the final monochromatic
$C_n^{2,c}$ we are constructing, so on these segments there will not
be any short diagonals guaranteed. However, since the paths $P_i^R$
are in $G^R$, their total length is indeed a constant depending on
$\eps$ only. We remove the internal vertices of the paths $P_i$ from
our graph; let us denote their total number by $c_1$. Thus on the
remaining segments we need exactly $n-c_1$ vertices.

Furthermore, we remove some more vertices from each $(V_1^i, V_2^i,
V_3^i), 1\leq i \leq l_1$ to achieve super-regularity in red in all
of the three pairs. Finally we remove some more vertices to get a
balanced super-regular tripartite graph. The number of remaining
vertices in the union of the clusters in the triangles is still
between $n$ and $(1+\eta)n$ using (\ref{para}). For simplicity we
still denote the clusters by $V_j^i$.

Finally we will lift the triangles back to almost square-paths in
the original graph using Lemma~\ref{blow-up}. Let us denote by
$(1-\eta')$ the ratio of $n$ and the number of remaining vertices in
the union of the clusters in the triangles, so $0\leq \eta'\ll 1$.
Let us use Lemma \ref{blow-up} in each balanced super-regular
tripartite graph $(V_1^i, V_2^i, V_3^i), 2\leq i \leq l_1$ with
$p_i=\lfloor (1-\eta')3|V_1^i|\rfloor$ to connect the two endpoints
of the connecting paths $P_{i-1}$ and $P_i$ with an almost
square-path of length $p_i$. Finally we use Lemma \ref{blow-up} one
more time in the balanced super-regular tripartite graph $(V_1^1,
V_2^1, V_3^1)$ with a $p_1$ value that makes the total length
exactly $n$, to connect the two endpoints of the connecting paths
$P_{\ell_1}$ and $P_1$ with an almost square-path of length $p_1$.
This is possible since this $p$ value is less than $\lfloor
(1-\eta')3|V_1^1|\rfloor$ only by a constant. Putting together the
almost square-paths within the triangles with the connecting paths
we get the red almost square-cycle of length $n$ with only a
constant number of short diagonals missing. \qed

\end{document}